\documentclass[12pt]{amsart}

\usepackage{amsmath, amssymb, graphics}

\usepackage{enumitem}

\newtheorem{theorem}{Theorem}[section]

\numberwithin{equation}{section}

\begin{document}

\author[Aimo Hinkkanen]{Aimo Hinkkanen}
\address{Department of Mathematics, University of Illinois Urbana-Champaign,
 Urbana, IL 61801 USA }
\email{aimo@illinois.edu}
\author[Joseph Miles]{Joseph Miles}
\address{Department of Mathematics, University of Illinois Urbana-Champaign,
 Urbana, IL 61801 USA}
\email{j-miles3@illinois.edu}

\title[Maximum and average valence]{Maximum and average valence of meromorphic functions}

\dedicatory{Dedicated to the memory of Lawrence Zalcman}

\subjclass[2010]{Primary 30D35}

\keywords{Meromorphic functions, Value distribution, Nevanlinna theory}

\maketitle

\begin{abstract}
If $f$ is a meromorphic function from the complex plane ${\mathbb C}$ to the extended complex plane $\overline{  {\mathbb C} }$, for $r > 0$ let $n(r)$ be the maximum number of solutions in $\{z\colon  |z| \leq r \}$ of $f(z) = a$ for $a \in \overline{  {\mathbb C} }$, and let $A(r,f)$ be the average number of such solutions.  Using a technique introduced by Toppila, we exhibit a meromorphic function for which $\liminf_{r\to\infty} n(r)/A(r,f) \geq 1.07328$.
\end{abstract}

\section{Introduction}

Let $f$ be a non-constant meromorphic function in the complex plane 
${\mathbb C}$. We denote the extended complex plane or Riemann sphere ${\mathbb C}\cup 
\{\infty\}$ by $\overline{  {\mathbb C} }$. We assume that the reader is familiar with Nevanlinna's value distribution theory as explained, for example, in \cite{H}. 

For $z\in {\mathbb C}$ and $r>0$, we write 
$B(z,r) = \{w\in    {\mathbb C} \colon |w-z| <  r \}$,
$\overline{B}(z,r) = \{w\in    {\mathbb C} \colon |w-z| \leq r \}$,
and 
$S(z,r) = \{w\in    {\mathbb C} \colon |w-z| = r \}$.  For convenience we set $B(z,0) = \{z\}$ (instead of $\emptyset$). 

If $r>0$ and $a\in \overline{  {\mathbb C} }$, we denote by $n(r,a)$ the number of points $z$ with $|z|\leq r$ and $f(z)=a$, with due count of multiplicity. We write
$$
n(r) = \sup \{n(r,a) \colon a\in  \overline{  {\mathbb C} }  \} 
$$
for the maximum valence of $f$ in the closed disk $\overline{B}(0,r)$. 
 The average valence of $f$ in $\overline{B}(0,r)$ with respect to the spherical metric in the image is given by
$$
A(r) = A(r.f) =  \frac{1}{\pi}  \int_0^{2\pi} \int_{0}^{r}  \frac{ | f'(te^{i\theta}) |^2  } {  ( 1 +  | f(te^{i\theta}) |^2  )^2  } \, t\, dt\, d\theta  .
$$

Since $n(r)\geq A(r)$, we have
$$
\liminf_{r\to\infty} \frac{ n(r) } { A(r) } \geq 1 .
$$
Hayman and Stewart \cite{HS} proved that
$$
\liminf_{r\to\infty} \frac{ n(r) } { A(r) } \leq e  .
$$
Miles \cite{Mil} obtained the improvement
$$
\liminf_{r\to\infty} \frac{ n(r) } { A(r) } \leq e - 10^{-28}  .
$$

Gary \cite{G}  has shown that for every finite subset $B$ of $ \overline{ {\mathbb C} }$, the quantity $n_B(r) = \max \{ n(r,a) \colon a\in B \}$ satisfies 
$\liminf_{r\to\infty} n_B(r)/A(r) \leq 2.65$. 

Toppila \cite{T} exhibited a meromorphic function with
$$
\liminf_{r\to\infty} \frac{ n(r) } { A(r) } \geq \frac{80}{79} \approx 1.01266 , 
$$
showing that this lower limit can be greater than $1$. Toppila's proof omits many details, so that for clarity alone, it may be worthwhile to give a more precise exposition of his method. We will do this, observing that parameters can be used in the construction, so that Toppila's example amounts to only one particular choice of the parameters. Using different choices, we obtain the following result.

\begin{theorem} \label{th1}
There exists a non-constant meromorphic function $f$ in the plane for which
\begin{equation} \label{1}
\liminf_{r\to\infty} \frac{ n(r) } { A(r) } \geq 1.07328 .
\end{equation} 
\end{theorem}

We note \cite{Mil1} that there exist meromorphic functions  for which \hfil\break 
$\limsup_{r\to\infty}   n(r) / A(r)  = \infty$.

\section{Preliminaries}

Recall that the chordal metric on $\overline{  {\mathbb C}  }$ is defined by
$$
k(z,a) = \frac{    |z-a| } {  \sqrt{ 1 + |z|^2 }     \sqrt{ 1 + |a|^2 }     } \leq \min\{  |z-a| ,1\} 
$$
for $z,a\in {\mathbb C}$ and
$$
k(z,\infty) = \frac{  1 } {  \sqrt{ 1 + |z|^2 }  }  .
$$
For $a\in   \overline{  {\mathbb C} }$ and $r>0$ 
we write
$$
D(a,r) = \{ z\in \overline{  {\mathbb C} } \colon k(z,a) < r \}  .
$$
Note that if $0<  \alpha  <1$ and $a\in  {\mathbb C} \setminus \{0\} $,  then
\begin{equation} \label{2}
B(a,  \alpha  |a|  ) \subset D(a, \alpha  ) 
\end{equation}
since if $|z-a|<  \alpha  |a|$, then
$$
k(z,a) < \frac{   \alpha  |a|     } {    \sqrt{ 1 + |z|^2 }     \sqrt{ 1 + |a|^2 }      } <  \alpha  .
$$
If $a=0$ and $0<  \alpha  <1$, then (\ref{2}) holds by our convention. 

Suppose that $0 <   \varepsilon   < 1$.  It is elementary that
$$
B(0,\sqrt{  \varepsilon/(1-\varepsilon)    } ) = D(0, \sqrt{ \varepsilon   } )
$$
has normalized spherical area $  \varepsilon $ in $\overline{  {\mathbb C} }$.  

If $L$ is a rotation of $\overline{  {\mathbb C} }$, then
$$
k(L(z),L(a)) = k(z,a)
$$
for all $z,a\in \overline{  {\mathbb C} }$, and so
\begin{equation} \label{3}
L(D( a, r     )) = D(L(a), r ) . 
\end{equation}
It follows that for every $a \in \overline{  {\mathbb C} }$, the disk $D(a,\sqrt{ \varepsilon   } )$ has area $\varepsilon$.  

For small $\varepsilon>0$, we now define three disks, centered at $0$, $\infty$, and $1$, each with 
area ~$\varepsilon$:
$$
\Omega_1 = D(0, \sqrt{ \varepsilon   } )  , \,\, 
\Omega_2 = D( \infty , \sqrt{ \varepsilon   }  )   , 
\,\, \text{ and } \,\,
\Omega_3 = D(1, \sqrt{ \varepsilon   }  )   . 
$$
We next define subdisks with the same centers but with half the chordal radius:
$$
\tilde{\Omega}_1 = D(0, \sqrt{ \varepsilon   }/2 )  , \,\, 
\tilde{\Omega}_2 = D( \infty , \sqrt{ \varepsilon   }/2  )   , 
\,\, \text{ and } \,\,
\tilde{\Omega}_3 = D(1, \sqrt{ \varepsilon   }/2  )   . 
$$
Note that if $L(z)=1/z$, a rotation of $\overline{  {\mathbb C} }$, then by 
(\ref{3}), we have
\begin{equation} \label{4}
L(  \Omega_1  ) = \Omega_2, \,\, 
L(  \Omega_2  ) = \Omega_1,  
\,\, \text{ and } \,\,
L(  \Omega_3  ) = \Omega_3   ,
\end{equation}
with the same relations holding for the $\tilde{\Omega}_j$. 
We set
$$
X = \overline{  {\mathbb C} }  \setminus (   {\Omega}_1 \cup  {\Omega}_2 \cup  {\Omega}_3       ) 
$$
and note that $m(X) = 1 - 3 \varepsilon$ where $m$ is the normalized spherical area measure on $ \overline{  {\mathbb C} }$. 

The M\"{o}bius transformation
$$
T(z) = \frac{  z + \frac{1}{3}  } {  1 + \frac{1}{3} z      } 
$$
plays a central role in our construction. Note that
$$
T^{-1}(z) = \frac{  z - \frac{1}{3}  } {  1 - \frac{1}{3} z      }   .
$$

For $a\in  {\mathbb C}  $ and $0 < r_1 < r_2$, let 
$$
A(a,r_1,r_2) =  \{ z\in {\mathbb C} \colon r_1 < |z-a| < r_2 \}
$$ 
be the annulus with center $a$, inner radius $r_1$, and outer radius $r_2$. 

Suppose that $\delta>0$ is small.  A direct calculation shows that
\begin{equation} \label{7}
T ( S ( 0, 1 + \delta ) ) \subset A( 0, 1 + \delta/3 , 1 + 3 \delta )  ,
\end{equation}
\begin{equation} \label{8}
T ( S ( 0, 1 - \delta ) ) \subset A( 0  , 1 - 3 \delta , 1 - \delta/3 )  ,
\end{equation}
\begin{equation} \label{9}
T^{-1} ( S ( 0, 1 + \delta ) ) \subset A( 0, 1 + \delta/3 , 1 + 3 \delta )  ,
\end{equation}
and
\begin{equation} \label{10}
T^{-1} ( S ( 0, 1 - \delta ) ) \subset A( 0  , 1 - 3 \delta , 1 - \delta/3 )  .
\end{equation}

For later use we record some properties of $T$ on the circle $S(0,1)$. We have
\begin{equation} \label{8a}
T \left( -\frac{1}{3} + i \frac{ 2 \sqrt{2} } { 3 }  \right) =
\frac{1}{3} + i \frac{ 2 \sqrt{2} } { 3 }  , 
\end{equation}
\begin{equation} \label{9a}
T \left( \frac{1}{3} + i \frac{ 2 \sqrt{2} } { 3 }  \right) =
\frac{7}{9} + i \frac{ 4 \sqrt{2} } { 9 }  , 
\end{equation}
\begin{equation} \label{10a}
0.60817 < 1 - \frac{  \arctan (   2 \sqrt{2}    )      } {  \pi } < 0.60818  , 
\end{equation}
and
\begin{equation} \label{11a}
0.78365  < 1 - \frac{  \arctan    \frac{ 4 \sqrt{2} } { 7 }      } {  \pi } < 0.78366  . 
\end{equation}

An elementary calculation shows that
$$
\frac{ {\rm Im}\,  T(  x + i \sqrt{ 1 - x^2 }     )     } {  {\rm Re}\,  T(  x + i \sqrt{ 1 - x^2 }     )       }  = \frac{  8 \sqrt{ 1 - x^2 }    } {   10 x + 6      } , \quad  -1 \leq x \leq 1, \,\, x\not= -\frac{3}{5}  .
$$

Set
$$
\alpha (x) = {\rm Arg} \, \,  T (   x + i \sqrt{ 1 - x^2 }     )   = 
\left\{ \begin{matrix}  \pi + \arctan   \frac{  8 \sqrt{ 1 - x^2 }    } {   10 x + 6      }  & , & -1\leq x < -\frac{3}{5} \\   \arctan   \frac{  8 \sqrt{ 1 - x^2 }    } {   10 x + 6      } & , & -\frac{3}{5} < x \leq 1     \\ \pi/2 & , & x = -\frac{3}{5}  \end{matrix}     \right.
$$
where ${\rm Arg}$ is the principal argument and by $\arctan$ we mean the principal inverse tangent taking values in $(-\pi/2,\pi/2)$. We have
$$
\alpha' (x) = \frac{ -8 (10+6x)     } {  (  (10 x + 6  )^2 + 64 (1 - x^2 )    ) \sqrt{ 1 - x^2 }      } , \quad -1<x<1  . 
$$
It is evident that $\alpha(x)$ decreases from $\pi$ to $0$ on $[-1,1]$ and a calculation shows that 
$|\alpha'(x)|$ is largest near $x=-1$. Thus for small $s >0$, if $x_1$ and $x_2$ are in $[-1,1]$ and $|x_1-x_2|= s $, then
\begin{equation} \label{12a}
| \alpha(x_1) - \alpha(x_2) | \leq  \alpha(-1) - \alpha(-1+ s )  < 4 \sqrt{ s }    . 
\end{equation}

\section{Construction of the Example}

\subsection{The function $f$} 
We now require $\delta$ to be so small that 
\begin{equation} \label{11}
(1 - 3 \delta)^2 -  (1 - \sqrt{ \varepsilon   } )^2   >   \sqrt{ \varepsilon   }  .
\end{equation}
This, of course, implies that
\begin{equation} \label{12}
1 - 3 \delta > 1 - \sqrt{ \varepsilon   }  , 
\quad \text{ or } \quad
 \delta <   \frac{ 1 } { 3 }  \sqrt{ \varepsilon   }   .
\end{equation}

For $C >1$ (we will later choose $C = 2.28228$), we let $k_0$ be a positive integer so that with 
$$
N_k = [ C^{k+k_0} ]
$$
for integers $k\geq 5$ we have
\begin{equation} \label{14} 
\sum_{k=5}^{\infty} \frac{1} {  (1+\delta)^{N_k}       } 
 < \frac{1}{6}   \sqrt{ \varepsilon   }  .  
\end{equation}
Here for real $x$ we denote by $[x]$ the greatest integer  $\leq x$. Note that this implies that
\begin{equation} \label{16a}
(1 - \delta)^{N_k} 
 < \frac{1}{6}   \sqrt{ \varepsilon   } , \quad  k\geq 5   .  
\end{equation}

For $k\geq 5$ let $D_k$ be the disk
$$
D_k = B(  k e^{i  \alpha_k  }  , 3/2        )  
$$
where the $\alpha_k $ are chosen so that the Euclidean distance between any two such disks is at least $1$. Write
$$
L_{k}(z) = \frac{2}{3} \left(    e^{  -  i  \alpha_k  } z - k \right)  ,
$$
so that
$$
L_k ( D_k) = B(0,1)  .
$$
Let 
$$
T_k(z) = T(L_k(z)) .
$$
Note that
$$
T_k \left(  \left( k-\frac{1}{2} \right) e^{i  \alpha_k }   \right) = T \left( L_k \left( \left( k-\frac{1}{2} \right) e^{i  \alpha_k  } \right) \right)  = T(-1/3) = 0 .
$$
We shall refer to $\left( k-\frac{1}{2} \right) e^{i  \alpha_k } $ as the {\bf special point} in $D_k$.

Define
$$
R_k(z) = (T_k(z))^{  N_k  }  . 
$$
Let
$$
S_k(z) = 1 - \frac{1} {   R_k(z)    }  = 1 - \frac{1} {  (T_k(z))^{  N_k  }     } .
$$
We set
\begin{equation} \label{17a}
f(z) = \frac{   \displaystyle \prod_{ k=5 \atop k\,\, \text{odd}    }^{\infty}   S_k(z)     } {   \displaystyle  \prod_{k=6 \atop  k\,\, \text{even}    }^{\infty}   S_k(z)        }  .
\end{equation}
Note that for odd $k$, the function $f$ has a pole of multiplicity $N_k$ at the special point $\left( k-\frac{1}{2} \right) e^{i  \alpha_k } $ in $D_k$ and has $N_k$ simple zeros on $\partial D_k$. 
The distribution of these zeros on the boundary of $D_k$ is not uniform, but rather is skewed due to the effect of the M\"{o}bius transformation $T$.
For even $k$, the function $f$ has a zero of multiplicity $N_k$ at the special point   in $D_k$ and has $N_k$ simple poles on $\partial D_k$, again not uniformly distributed. 

\subsection{A sketch of the proof} 
Before presenting the rather intricate details of the argument establishing the required property of $f$, we give a brief overview of the strategy behind the example.  Consider an odd $k$.  Our function has a pole of multiplicity $N_k$ on $S ( 0, k- 1/2)$.  Importantly (see (\ref{50a})), for most values of $r$ in $( k - 1/2 ,  k + 3/2)$ and for most values of $a$ in $\overline{{\mathbb C} }$, the number of solutions of $f(z) = a$ inside   $S (0 , r)$ and near $D_k$ is very close to the number of solutions of $f(z) = 0$ inside $S(0,r)$ and on the boundary of $D_k$.  Thus for $k - 1/2 < r < k + 3/2$, the contribution to $A(r,f)$ from $z$ close to $D_k$ lags behind the number of poles of $f$ in $D_k$ and only ``catches up'' at $r = k + 3/2$.  However on   $S( 0, k + 1/2)  = S(0, k+ 1 - 1/2)$ lies the special point in $D_{k +1}$, at which $f$ has a zero of multiplicity $N_{k + 1}$.  The argument is then repeated (see (\ref{51a})) with the roles of the zeros and poles reversed. Thus $A(r,f)$  never ``catches up'' to $n(r)$.  We now turn to the details of the example.

\subsection{Certain circles and annuli} 
We  define a circle $C_a^k$ just outside $\partial D_k$ and a circle $C_b^k$ just inside $\partial D_k$ that play an important role in our construction.  The circle $C_a^k$ is defined by
$$
L_k ( C_a^k ) = T^{-1} (  S(0, 1 + \delta    )      ) ,
$$
so that 
$$
T_k ( C_a^k ) = T(   L_k ( C_a^k )     )  =   S(0, 1 + \delta    )   .
$$
The circle $C_b^k$ is defined by
$$
L_k ( C_b^k ) = T^{-1} (  S(0, 1 - \delta    )      ) ,
$$
so that
$$
T_k ( C_b^k ) =    S(0, 1 - \delta    )   .  
$$
Note that the centers of $C_a^k$ and $C_b^k$ have argument $\alpha_k$. 

The relations  (\ref{9}) and  (\ref{10}) imply that 
\begin{equation} \label{15}
L_k ( C_a^k ) \subset A(0, 1 + \delta /3 , 1 + 3 \delta ) 
\end{equation}
and
\begin{equation} \label{16}
L_k ( C_b^k ) \subset A(0 , 1 - 3 \delta , 1 - \delta /3  ) .  
\end{equation}

Let $G_k$ be the disk inside $C_b^k$. Let $H_k$ be the disk inside $C_a^k$. Let $P_k$ be the region of $  {\mathbb C} $  outside $C_b^k$. Let $J_k$ be the region of  $ {\mathbb C} $  outside $C_a^k$. Set
$$
U_k = H_k \cap P_k .  
$$
Note that $U_k$ is  bounded by the two  circles $C_a^k$ and $C_b^k$. 
By  (\ref{15}) and  (\ref{16})  we have
\begin{equation} \label{20a}
U_k   \subset A \left(  k e^{i  \alpha_k  } , \frac{3}{2} -  \frac{9}{2} \delta , \frac{3}{2} +  \frac{9}{2}  \delta \right)  .  
\end{equation}
Also
$$
D_k \subset \bigcap_{j\not= k}   J_j  .
$$

Certainly, with the notation $E^k=\{ z^k \colon z\in E\}$, 
$$
R_k ( C_b^k ) 
= ( T_k(C_b^k) )^{N_k} 
= S (0, (1 - \delta )^{N_k} )  .
$$
For $z\in C_b^k$,  by (\ref{16a}) we have
$$
k( R_k(z), 0) < | R_k(z)    | =  (1 - \delta )^{N_k}   <  \sqrt{ \varepsilon   } /6 ,
$$
implying that 
\begin{equation} \label{17}
R_k ( C_b^k ) \subset  D\left(  0,    \sqrt{ \varepsilon   } /6  \right)    \subset   \tilde{\Omega}_1  .  
\end{equation}
Similarly, 
\begin{equation} \label{18}
R_k ( G_k )  =   B(0,     (1 - \delta )^{N_k}      ) \subset   D\left(  0,    \sqrt{ \varepsilon   } /6  \right)    \subset       \tilde{\Omega}_1  .  
\end{equation}

From (\ref{14}) we have 
\begin{equation} \label{19}
R_k ( C_a^k )  =   S(0,     (1 + \delta )^{N_k}      )    \subset  D\left(  \infty ,    \sqrt{ \varepsilon   } /6  \right)      \subset \tilde{\Omega}_2    
\end{equation}
and 
\begin{equation} \label{20}
R_k ( \overline{ J_k } )  =  \{ w \in \overline{  {\mathbb C} } \colon |w| \geq  (1 + \delta )^{N_k}   \}     \subset  D\left(  \infty ,    \sqrt{ \varepsilon   } /6  \right)   \subset \tilde{\Omega}_2  .  
\end{equation}
We take $\overline{ J_k } $  to be the closure of $J_k$ in ${\mathbb C}$.

We conclude that 
\begin{equation} \label{21}
S_k ( C_b^k ) =  S( 1, (1 - \delta )^{-N_k}          )       \subset  D\left(  \infty ,    \sqrt{ \varepsilon   } /3  \right)           \subset \tilde{\Omega}_2    
\end{equation} 
and 
\begin{equation} \label{22}
S_k ( \overline{ G_k } )  =   \{ w\in \overline{  {\mathbb C} }  \colon |w-1| \geq      (1 - \delta )^{-N_k}      )         \subset  D\left(  \infty ,    \sqrt{ \varepsilon   } /3  \right)         \subset \tilde{\Omega}_2    . 
\end{equation}

We also have
\begin{equation} \label{23}
S_k ( C_a^k )  =  S(  1,   (1 + \delta )^{-N_k}       )      \subset  D\left(  1 ,    \sqrt{ \varepsilon   } /6  \right)       \subset \tilde{\Omega}_3    
\end{equation}
and
\begin{equation} \label{24}
S_k (  \overline{ J_k } )  =  \overline{B}(  1,   (1 + \delta )^{-N_k}       )    \subset  D\left(  1 ,    \sqrt{ \varepsilon   } /6  \right)         \subset \tilde{\Omega}_3    . 
\end{equation}
Note that (\ref{14}) and (\ref{24})  guarantee the convergence of the infinite products appearing in (\ref{17a}).

Suppose that $z\in \cap_{j=5}^{\infty} \overline{ J_j } $. Denoting the principal logarithm by ${\rm Log}$, we have
$$
{\rm Log}\, f(z) =  \sum_{\scriptstyle j \,\, {\rm odd} \atop \scriptstyle j\geq 5}  {\rm Log}\, S_j(z) 
-  \sum_{\scriptstyle j \,\, {\rm even} \atop \scriptstyle j\geq 6}  {\rm Log}\, S_j(z)  
$$
and from (\ref{24}) 
$$
\left|     {\rm Log}\, S_j(z)        \right|    \leq 2 \left| S_j(z) - 1       \right| < \frac{2}{    (1 + \delta )^{N_j}    }  .
$$
From (\ref{14})  we conclude that  
$$
|  {\rm Log}\, f(z) |  
< \frac{1}{3} \sqrt{ \varepsilon}  ,  
$$
implying that
\begin{equation} \label{25a} 
 | f(z) - 1|< \frac{1}{2} \sqrt{ \varepsilon} . 
\end{equation}
Thus 
$$
k(f(z),1) \leq |f(z)-1| < \frac{1}{2} \sqrt{ \varepsilon} ,
$$
or  
\begin{equation} \label{30}
f(z) \in \tilde{\Omega}_3  \quad \text{ for all } \,\, z\in \bigcap_{j=5}^{\infty}  \overline{J_j}  . 
\end{equation}
In particular, since $\overline{U_k} \subset \cap_{\scriptstyle j=5 \atop \scriptstyle  j\not= k}^{\infty} \overline{ J_j }$, we have
\begin{equation} \label{27a}
f(  C_a^k   )  \subset f\left( \bigcap_{j=5}^{\infty}   \overline{J_j}    \right)    \subset \tilde{\Omega}_3  . 
\end{equation}

If $k$ is odd, then in the sense of the chordal metric, $f$ closely approximates $S_k$ on $\tilde{D}_k := B(   k e^{i  \alpha_k  }   , 2)$.  Specifically, consider $z$ in $\tilde{D}_k $ with 
$S_k(z) \not=   \infty$. (Of course $S_k(z) = 0$ or $\infty$ in $\tilde{D}_k $ 
 if, and only if, $f(z) = 0$ or~$\infty$.) 
  Note that for $z\in \tilde{D}_k $ we have $z \in \cap_{j \not= k } \overline{ J_j } $. 
Writing 
$$
f(z) = S_k(z) g_k(z) , 
$$ 
we may apply the discussion leading to
(\ref{25a}) to $g_k$ to conclude that 
$$
| g_k(z) - 1 | < \frac{1}{2} \sqrt{\varepsilon}  . 
$$
Since
$$
|f(z)-S_k(z)| = |S_k(z)| |g_k(z) - 1|  ,
$$
using (\ref{2}) we have
\begin{eqnarray} \label{33}
f(z) \in B( S_k(z) , |S_k(z)| |g_k(z) - 1| ) & \subset &  B \left( S_k(z) , \frac{1}{2} \sqrt{\varepsilon} |S_k(z)|  \right) \notag
\\
& \subset &  D \left( S_k(z) , \frac{1}{2} \sqrt{ \varepsilon}  \right) ,  
\end{eqnarray}
for all $z\in \tilde{D}_k$ with $f(z) \not= \infty$. 
The fact that (\ref{33}) holds for all $z\in \overline{ U_k }$, including those $z$ for which $S_k(z)=0$, plays a key role in our comparison of the solutions of $S_k(z)=a$ to the solutions of $f(z)=a$. 

If $z \in  \overline{ G_k }$ and $f(z) \not= \infty$, i.e., $z$ is not the special point in $D_k$, applying (\ref{22}) and (\ref{33}), we have
$$
k(\infty , f(z) ) \leq k(\infty , S_k(z) ) + k( S_k(z) , f(z) ) < \frac{1}{3}  \sqrt{ \varepsilon} + \frac{1}{2}  \sqrt{ \varepsilon} < 
\sqrt{ \varepsilon} .
$$
Clearly, this inequality holds if $f(z) = \infty$  as well. 
We conclude that 
\begin{equation} \label{29a}
f(  C_b^k   )  \subset f(    \overline{G_k}    )    \subset {\Omega}_2   
\end{equation}
for odd values of $k$. 

If $j$ is even, we note that $S_j(z)$ is a factor in the numerator of $1/f$. The above argument shows that 
$$
\frac{1}{f} (     \overline{G_j }     ) \subset  {\Omega}_2  , 
$$
or
\begin{equation} \label{30a}
f(  C_b^j   )  \subset f(    \overline{G_j}    )    \subset {\Omega}_1  .     
\end{equation}

For $k\geq 5$ 
let $n_k(a,f)$ denote the number of solutions of $f(z)=a$ in $U_k$ counting multiplicities. 
For $g=S_k$ or $g=f$, let $n_k(r,a,g)$ be the number of solutions of $g(z)=a$ in $U_k \cap \overline{B}(0,r)$. 

For odd $k$, we know from (\ref{29a}) that 
$f(C_b^k) \subset \Omega_2$ and from (\ref{27a}) that 
$f(C_a^k) \subset \Omega_3$.  By the Argument Principle, $n_k(a,f)$ is constant on 
${\overline {\mathbb C} } \setminus ( \Omega_2 \cup \Omega_3 ) \supset X $. Since 
$0\in {\overline {\mathbb C} } \setminus ( \Omega_2 \cup \Omega_3 )$ and 
$n_k(0,f) = N_k $, we conclude for odd $k$ that 
\begin{equation} \label{31a}
n_k(a,f)=N_k
\end{equation} 
for all $a\in {\overline {\mathbb C} } \setminus ( \Omega_2 \cup \Omega_3 ) \supset X $.  

If $j$ is even, we have from (\ref{27a}) that $f(C_a^j) \subset \Omega_3$ and from (\ref{30a}) that $f(C_b^j) \subset \Omega_1$. Again we conclude that 
\begin{equation} \label{32a}
n_j (a,f) = n_j (\infty, f) = N_j 
\end{equation}
for all $a\in {\overline {\mathbb C} } \setminus ( \Omega_1 \cup \Omega_3 ) \supset X $.

We have 
$$
f\left( \bigcap_{k=5}^{\infty} \overline{ J_k } \right) \subset \tilde{\Omega}_3 
$$
by (\ref{30}). For odd $k$,
$$
f( \overline{ G_k } ) \subset {\Omega}_2 
$$
by (\ref{29a}). For even $k$, 
$$
f( \overline{ G_k } ) \subset {\Omega}_1 
$$
by (\ref{30a}). Since
$$
{\mathbb C} = \left( \bigcap_{k=5}^{\infty}  \overline{J_k}   \right) \cup \left( \bigcup_{k=5}^{\infty}  U_k  \right) \cup \left( \bigcup_{k=5}^{\infty}   \overline{G_k}  \right)  ,  
$$
we conclude that
\begin{equation} \label{33a}
f^{-1} (X) \subset \bigcup_{k=5}^{\infty} U_k  .  
\end{equation}

\subsection{Certain image areas} 
We may write $A(r,f)$ in the form 
$$
A(r,f) = \int_{  \overline {\mathbb C}     }  n(r,a,f) \, dm(a)  . 
$$
For a Lebesgue measurable set $Y\subset  \overline {\mathbb C}  $, let
$$
A^Y (r,f) = \int_{ Y }  n(r,a,f) \, dm(a)  
$$
and
$$
A_k^Y (r,f) = \int_{ Y }  n_k(r,a,f) \, dm(a)  .  
$$
We write $A_k (r,f) $ for $A_k^Y (r,f) $ if $Y=  \overline {\mathbb C}  $.  

In view of (\ref{33a}), if $Y\subset X$ is measurable, then
\begin{equation} \label{34a}
A^Y (r,f) = \sum_{ k=5 }^{\infty}    A_k^Y (r,f) .  
\end{equation}
If $U_k \subset  \overline{B}(0,r) $, then from (\ref{31a}) and (\ref{32a})  we have   
$$
A_k^Y(r,f) = N_k\,  m(Y)  
$$
for measurable $Y\subset X$. 
If $U_k \cap \overline{B}(0,r) = \emptyset$, evidently $A_k^Y(r,f)=0$ for every measurable $Y$. Note that for all such $Y$ and $r$ we have
\begin{equation} \label{35a}
A(r,f) \leq  A^Y (r,f) + n(r) \, (1 - m(Y))  .  
\end{equation}

For measurable $Y\subset X$, we turn our attention to $A_k^Y (r,f) $ for $r$  such that 
$$
\emptyset \not=  \overline{B}(0,r)  \cap U_k \subsetneq U_k . 
$$

\subsection{Intersections of certain circles} 
For $k\geq 5$ and real $\beta$ with 
\begin{equation} \label{36a}
| \beta | \leq \frac{3}{2}  - \frac{3}{2}  \sqrt{ \varepsilon }  ,  
\end{equation}
we consider the circle
$S(0, k+\beta)$. For $\beta>0$, by (\ref{12}) we have
$$
k+\beta \leq k+ \frac{3}{2} -  \frac{3}{2} \sqrt{ \varepsilon } < k + \frac{3}{2} - \frac{9}{2} \delta 
$$
and for $\beta<0$, we have
$$
k+\beta \geq k - \frac{3}{2} + \frac{3}{2} \sqrt{ \varepsilon } >  k - \frac{3}{2} + \frac{9}{2} \delta . 
$$
Let $C'$ be the circle $S \left(k e^{i \alpha_k} , \frac{3}{2} + \frac{9}{2} \delta \right)$ and $C''$ the circle $S \left(k  e^{i \alpha_k} , \frac{3}{2} - \frac{9}{2} \delta  \right)$. 
Note that by (\ref{20a}), $C_a^k$ lies in the annulus between $\partial D_k$ and $C'$; likewise $C_b^k$ lies in the annulus between $\partial D_k$ and $C''$. 
By our choice of $\beta$, the circle $S(0,k+\beta )$  intersects $C''$ (as well as the other four circles $C_b^k$, $\partial D_k$, $C_a^k$, and $C'$) in two points.

For $1\leq j\leq 5$, let
$$
z_j = (k+\beta) e^{ i ( \alpha_k +  \theta_j ) }  ,  
$$
where
$$
0 < \theta_1 < \theta_2 < \theta_3 < \theta_4 < \theta_5 < \frac{\pi}{2}  , 
$$
be  points where $S(0,k+\beta )$  intersects
$C''$, $C_b^k$, $\partial D_k$, $C_a^k$, and $C'$, respectively.  Let $\gamma_j$ be the closed arc of 
$S(0,k+\beta )$  from $z_j$ to $z_{j+1}$, for $1\leq j\leq 4$. 

Note that the points 
$\tilde{z}_j = (k+\beta) e^{ i ( \alpha_k -  \theta_j ) } ,$
for $1\leq j\leq 5$, are the other points of intersection of $S(0,k+\beta )$ with
$C''$, $C_b^k$, $\partial D_k$, $C_a^k$, and $C'$, respectively, since the centers of $C_b^k$ and $C_a^k$ have argument $\alpha_k$. For $1\leq j\leq 4$, let $\tilde{\gamma}_j$  be the closed arc of $S(0,k+\beta )$   from    $\tilde{z}_{j+1}$ to $\tilde{z}_j$.   

Let $\gamma_a$ be the arc of $C_a^k$ from $z_4$ to $\tilde{z}_4$ that intersects the ray from $0$ to $k  e^{i \alpha_k}$. Let $\gamma_b$ be the arc of $C_b^k$ from  $\tilde{z}_2$ to $z_2$   that intersects the ray from $0$ to $k  e^{i \alpha_k}$. 
Let $\gamma$ be the closed path consisting of $\gamma_2 \cup \gamma_3$, $\gamma_a$, $\tilde{\gamma}_3 \cup \tilde{\gamma}_2$, and $\gamma_b$, oriented so that the interior $V$ of $\gamma$ lies on the left.  Note that
$$
V = U_k \cap B(0,k+\beta )  .
$$
We wish to examine the mapping  $S_k$ on $V$. Let 
$
\Gamma = S_k(\gamma) , \, \Gamma' = S_k(\gamma_2 \cup \gamma_3) ,  \, \Gamma'' = S_k(\tilde{\gamma}_3 \cup \tilde{\gamma}_2)   , \, \Gamma''' = S_k(\gamma_a)   , 
\, \text{ and } \,  \Gamma'''' =  S_k(\gamma_b)   .  
$ 

Let $A_k^*$ be the annulus
$$
A_k^* = A( 1, (1+\delta)^{-N_k} ,  (1-\delta)^{-N_k}       )  .
$$
Note that $0\in A_k^* $ and
$$
A_k^*  \supset   \overline{{\mathbb C}} \setminus (    \Omega_2 \cup \Omega_3         )       \supset  X  .
$$

By (\ref{23}), $\Gamma'''$ lies on the inner boundary of $A_k^*$.  
By (\ref{21}), $\Gamma''''$ lies on the outer boundary of $A_k^*$. 
As $z$ moves along $\gamma_2 \cup \gamma_3$, the point $S_k(z)$ moves along $\Gamma'$ from the outer boundary of $A_k^*$  to its inner boundary. As $z$ moves along $\tilde{\gamma}_3 \cup \tilde{\gamma}_2$,  the point $S_k(z)$ moves along $\Gamma''$ from the inner boundary of $A_k^*$  to its outer boundary.

Suppose that 
$$
\frac{1}{ ( 1+ \delta)^{N_k}  }     < t_1<t_2 < \frac{1}{( 1- \delta)^{N_k}  }. 
$$
Then 
$$
C_{t_j} := \{ z \colon | S_k(z) - 1 | = t_j \}
$$
is a circle that lies between $C_a^k$ and $C_b^k$, i.e., lies in $U_k$, and whose center has argument $\alpha_k$. 
Furthermore $C_{t_2}$ lies inside $C_{t_1}$. Thus 
$| S_k(z) - 1 |$ is strictly decreasing on $\Gamma'$ and strictly increasing on $\Gamma''$. It follows that if $\tilde{C}$ is any circle with center $1$ lying in $A_k^*$, the arcs $\Gamma'$ and $\Gamma''$ each intersect $\tilde{C}$ in exactly one point. Thus if $a$ and $b$ both belong to $ \tilde{C} \setminus (   \Gamma' \cup \Gamma''         )    $, then
\begin{equation} \label{37a}
| {\rm Ind}_{\Gamma} (a) -  {\rm Ind}_{\Gamma} (b)           | \leq 1 , 
\end{equation}
where ${\rm Ind}_{\Gamma} (a)$ denotes the winding number of $\Gamma$ about $a$. 
Since $V = U_k \cap B(0,k+\beta )$, the Argument Principle gives
\begin{equation} \label{38a}
| n_k( k+\beta , a, S_k ) -  n_k( k+\beta , b, S_k )   | \leq 1    
\end{equation}
for all $a$ and $b$ in  $A_k^*$  with  $|a-1|=|b-1|$.  

We turn our attention to $\Gamma' \cap \Gamma''$.  It is evident that 
$S_k(z)$ is real if $z$ lies on the line through $0$ and $e^{i \alpha_k}$. For
$0<\theta < \pi$, we have
$$
S_k (  (k+\beta)  e^{i (\alpha_k + \theta ) }      ) = \overline{   S_k (  (k+\beta)  e^{i (\alpha_k - \theta ) }      )        } . 
$$
Thus  $\Gamma' $ and $ \Gamma''$ are complex conjugates of one another. Since each of these curves intersects a circle $\tilde{C} = S(1,t) \subset A_k^*$ exactly once, all points of intersection of $\Gamma' $ and $ \Gamma''$ in $A_k^*$ must lie on the real axis ${\mathbb R}$. 

We wish to obtain an upper estimate for the number of points of intersection of $\Gamma' $ and $ \Gamma''$ in $A_k^*$.
We begin by considering a continuous $\arg (w-1)$ for $w\in \Gamma'$. Recall that $\Gamma' = \{ S_k ( (k+\beta ) e^{ i ( \alpha_k + \theta ) } )  \colon \theta_2 \leq \theta \leq \theta_4 \}$. 
From the definition of $S_k$, we have
\begin{equation} \label{38b}
\arg (S_k(z) - 1 ) = \pi - N_k \, \arg T_k(z) = \pi - N_k \, \arg (T(L_k(z))) .
\end{equation}
Thus it is sufficient to consider a continuous $\arg (T(L_k(z)))$ for $z\in \gamma_2 \cup \gamma_3$. 

Let $C^* = T(L_k(S(0,k+\beta )))$. We note that
$
L_k(  (k+\beta ) e^{i \alpha_k }     ) = \frac{2}{3} \beta
$
and
$
L_k( - (k+\beta ) e^{i \alpha_k }     ) = - \frac{4}{3} k - \frac{2}{3} \beta .
$
A direct calculation yields
$$
T \left( \frac{2}{3} \beta \right) = 3\left(  1 - \frac{8} {  2\beta + 9  }      \right) \in (-1,1) 
$$
and
$$
T \left( - \frac{4}{3} k - \frac{2}{3} \beta \right) = 3\left(  1 + \frac{8} { 4k +  2\beta - 9  }      \right) \approx 3  .   
$$
Thus, for large $k$, the set $C^*$ is a circle with center $c$ in ${\mathbb R}$ where
$$
c \approx 3 - \frac{  12  } {  2\beta + 9    } \in (1,2) 
$$
and 
radius
$$
\rho \approx \frac{  12  } {  2\beta + 9    }  \in (1,2)  .
$$
Since $T$ maps the upper half plane $H$ to itself, we see that as $z$ traverses $S(0,k+\beta)$ in the counterclockwise direction, $T(L_k(z))$ traverses $C^*$ in the clockwise direction. We consider three cases.

If $\beta < -1/2$ and $\beta$ satisfies (\ref{36a}), then 
$T \left( \frac{2}{3} \beta \right) \in (-1,0)$, the circle $C^*$ surrounds $0$, and  
$\arg T(L_k( (k+\beta ) e^{i ( \alpha_k + \theta ) } ))$ is decreasing on $(0,\pi)$. 
Thus $\arg (S_k(  (k+\beta)  e^{i (\alpha_k + \theta ) } ) - 1 )$ is increasing on $(0,\pi)$ and in particular on $[ \theta_2,\theta_4 ]$. It follows that at every point $P$ of $\Gamma' \cap {\mathbb R} \cap A_k^*$, the arc $\Gamma'$ crosses ${\mathbb R}$  (either from the upper half plane to the lower half plane, or vice versa). The complex conjugate curve $\Gamma''$ also has this property. 

If $\beta = -1/2$, then $C^*$ contains $0$ and lies in the closed right half plane. As $z$ traverses $S(0,k+\beta)$ in the counterclockwise direction, $T(L_k(z))$ traverses $C^*$ in the clockwise direction. Thus $\arg (S_k(  (k+\beta)  e^{i (\alpha_k + \theta ) } ) - 1 )$ is increasing on $(0,\pi)$ and, as before, $\Gamma'$ crosses ${\mathbb R}$ at each point  of $\Gamma' \cap {\mathbb R} \cap A_k^*$, as does $\Gamma''$. 

If $\beta > -1/2$ and $\beta$ satisfies (\ref{36a}), $T \left( \frac{2}{3} \beta \right) \in (0,1) $  and $C^*$ lies in the open right half plane. Thus there exists $\theta^*$ in $(0,\pi)$ such that $\arg (T(L_k(  (k+\beta)  e^{i (\alpha_k + \theta ) }  )))$  is increasing on $(0,\theta^*)$ and decreasing on $(\theta^*,\pi)$. Evidently $
\arg (S_k(  (k+\beta)  e^{i (\alpha_k + \theta ) } ) - 1 ) 
$ 
is decreasing on $(0,\theta^*)$ and increasing on $(\theta^*,\pi)$. Thus 
 in the (quite unlikely) case that $\theta_2 < \theta^* < \theta_4$ and 
$
S_k(  (k+\beta)  e^{i (\alpha_k + \theta^* ) } ) - 1  
$ 
is purely real, there is a point $P$ (but only one such point) in $\Gamma' \cap {\mathbb R} \cap A_k^*$ at which $\Gamma'$ does not cross ${\mathbb R}$ but rather  is tangent to ${\mathbb R}$. Otherwise, $\Gamma'$ crosses ${\mathbb R}$ 
at each point of $\Gamma' \cap {\mathbb R}\cap A_k^*$. 

Consider the intervals 
$I^* = ( 1 - (1-\delta)^{-N_k} , 1 - (1+\delta)^{-N_k} )$
and
$I^{**} = ( 1 + (1+\delta)^{-N_k} , 1 + (1-\delta)^{-N_k} )$.
By the above observations for at least one of these intervals $I$, at all points $P$ in $I \cap \Gamma'  = I \cap \Gamma''$, the curve $\Gamma'$ and its conjugate curve $\Gamma''$ both cross $I$. For definiteness, suppose that $I^{**}$ has this property. Let $I^{**} \cap \Gamma'  =  I^{**}  \cap \Gamma'' = \{P_1,P_2,\dots ,P_m\}$ where $P_1<P_2<\cdots <P_m$. Let $P_0 = 1 + (1+\delta)^{-N_k} $ and $P_{m+1} = 1 + (1-\delta)^{-N_k} $. Let $Q_j = (P_j,P_{j+1})$, for $0\leq j\leq m$. Select $q_j$ in $Q_j$, for $0\leq j\leq m$. Because $\Gamma'$ and $\Gamma''$ both cross ${\mathbb R}$ at each $P_j$ for $1\leq j\leq m$, we have
\begin{equation} \label{39a}
| {\rm Ind}_{\Gamma} ( q_{j+1} ) -  {\rm Ind}_{\Gamma} ( q_j )           |   = 2 ,
\quad 0\leq j\leq m - 1 .
\end{equation}

We next obtain an upper bound on the total variation of a continuous $\arg (w-1)$ for $w\in \Gamma'$. By (\ref{38b}) it is sufficient to find an upper bound for the total variation of $\arg (T(L_k(z)))$ for $z\in \gamma_2 \cup \gamma_3$. Recalling that $T$ maps $H$ onto itself, we see that
$$
T(L_k(  \gamma_2 \cup \gamma_3 )) = C^* \cap \overline{ A(0,1-\delta,1+\delta)} \cap H  .
$$

Let $r_j e^{ i \varphi_j } $ with $|r_j-1| \leq \delta$ and $0<\varphi_j<\pi$ for $j=1,2$ 
be two arbitrary points in $T(L_k(  \gamma_2 \cup \gamma_3  ))$. By the Law of Cosines we have
$$
r_j^2 + c^2 - 2 c r_j \cos \varphi_j = \rho^2  , \quad j=1,2  .
$$
Subtracting we obtain
$$
r_2^2 - r_1^2 = 2c (r_2 \cos \varphi_2 - r_1 \cos \varphi_1 ) . 
$$
Rearranging, we obtain
$$
\frac{  (r_2 - r_1) (r_2 + r_1 )    } { 2 c } = \cos \varphi_2 -  \cos \varphi_1 
+ (r_2 - 1) \cos \varphi_2 - (r_1 - 1) \cos \varphi_1  .  
$$
Further rearrangement yields
$$
2 \sin \frac{  \varphi_1 + \varphi_2    } {2 }  \sin  \frac{ \varphi_1 - \varphi_2  } {2 } 
= \frac{  (r_2 - r_1) (r_2 + r_1 )    } { 2 c  } + (r_1 - 1) \cos \varphi_1 - (r_2 - 1) \cos \varphi_2  .  
$$
Thus
$$
\left|   \sin \frac{  \varphi_1 + \varphi_2    } {2 }  \sin  \frac{ \varphi_1 - \varphi_2  } {2 }         \right| \leq 
\frac{  2 \delta (  2 + 2\delta    )     } {  4 c } + \frac{\delta}{2} + \frac{\delta}{2} < 3\delta  .  
$$

Note that $\frac{ \varphi_1 + \varphi_2  } {2 }   \in (0,\pi)$ and $ \frac{ \varphi_1 - \varphi_2  } {2 }  \in (-\pi/2,\pi/2)$. 
We first suppose that
$\left|    \sin \frac{  \varphi_1 + \varphi_2    } {2 }    \right| > \sqrt{\delta}$. Then
$$
\frac{   |    \varphi_1 - \varphi_2   | } { \pi } \leq \left|   \sin  \frac{ \varphi_1 - \varphi_2  } {2 }         \right| < 3  \sqrt{\delta}  ,  
$$
implying that $   |    \varphi_1 - \varphi_2   | <  3 \pi  \sqrt{\delta} $. 

We next suppose that $\left|    \sin \frac{ \varphi_1 + \varphi_2  } {2 }    \right| \leq  \sqrt{\delta}$. We have either 
$\frac{  \varphi_1 + \varphi_2   } {2 } \leq \arcsin \sqrt{\delta} < 2 \sqrt{\delta} $ or
$\frac{  \varphi_1 + \varphi_2   } {2 } \geq \pi - \arcsin \sqrt{\delta}$. In the former case
we have 
$  |    \varphi_1 - \varphi_2   |  \leq  \varphi_1 + \varphi_2 < 4   \sqrt{\delta} $. In the latter case we have
$  \varphi_1 + \varphi_2 \geq 2\pi - 2 \arcsin \sqrt{\delta} $, implying that
$\varphi_j \geq \pi - 2 \arcsin \sqrt{\delta}$ for $j=1,2$ and 
$ |   \varphi_1 - \varphi_2     | \leq 2 \arcsin \sqrt{\delta} < 4 \sqrt{\delta}  $. 

Since $ \varphi_1 $ and  $\varphi_2 $ are arguments of arbitrary points in 
$T(L_k(  \gamma_2 \cup \gamma_3 ))$ and $\arg (T(L_k(z)))$ is either monotone on $\gamma_2 \cup \gamma_3$ or has one arc of increase and one arc of decrease, we conclude that the total variation of $\arg (T(L_k(z)))$  on $\gamma_2 \cup \gamma_3$ is at most $6\pi \sqrt{\delta}$. In view of (\ref{38b}), the total variation of $\arg (w-1)$ on $\Gamma'$ (equivalently of $\arg (S_k(z)-1)$ on $\gamma_2 \cup \gamma_3$) is at most $6\pi N_k \sqrt{\delta} < 2 \pi \sqrt{3} N_k \varepsilon^{1/4}$.

We now obtain an upper bound on $m$, the number of points $P_j$ in $\Gamma' \cap I^{**}$. Since $\arg (w-1)$ is monotone on $\Gamma'$ or has one arc of increase and one arc of decrease, as $w$ traverses $\Gamma'$ from $P_j$ to $P_{j+1}$, $1\leq j\leq m-1$, with at most one exceptional $j$, the change in a continuous $\arg (w-1)$ is $2\pi$. Thus the total variation of $\arg (w-1)$ on $\Gamma'$ is at least $(m-2)2\pi$, yielding
$$
(m-2)2\pi < 2\pi \sqrt{3} N_k \varepsilon^{1/4} .
$$

From (\ref{39a}) it follows that if $p$ and $q$ are in $I^{**} \setminus \Gamma$, then
$$
| {\rm Ind}_{\Gamma} (p) -  {\rm Ind}_{\Gamma} (q)           | \leq 2m < 2 \sqrt{3} N_k \varepsilon^{1/4}  + 4 .
$$
From (\ref{37a}) we conclude that if $a$ and $b$ are in $A_k^* \setminus \Gamma$, then
\begin{equation} \label{41b}
| {\rm Ind}_{\Gamma} (a) -  {\rm Ind}_{\Gamma} (b)           |   <
\eta(\varepsilon, k)      
\end{equation}
where to simplify notation we have set $\eta(\varepsilon, k) =4 N_k \varepsilon^{1/4} $.  

In the above discussion we have implicitly assumed that $m\geq 3$. If $m=1$ or $m=2$, then from (\ref{37a}) and (\ref{39a}) we conclude that the left hand side of (\ref{41b}) is bounded above by $2$ or $4$, respectively. If $m=0$, then ${\rm Ind}_{\Gamma} (a)$ is constant on $I^{**}$ and from (\ref{37a}) we conclude that the left hand side of (\ref{41b}) is bounded above by $1$.

By the Argument Principle we have from (\ref{41b}) that 
\begin{equation} \label{42a}
| n_k ( k+\beta , a,S_k ) - n_k ( k+\beta , b ,S_k )  | \leq \eta(\varepsilon, k)   
\end{equation}
for all $a$ and $b$ in $A_k^* $. In particular with $r = k+\beta $ we have
\begin{equation} \label{43a}
 n_k ( r , 0 , S_k ) - \eta(\varepsilon, k)  \leq n_k(r,a,S_k) \leq  n_k ( r , 0 , S_k ) + \eta(\varepsilon, k)    
\end{equation}
for all $a$ in $A_k^* $. 

From (\ref{43a})  we have 
\begin{eqnarray*}
&{}&
( n_k ( k+\beta , 0 , S_k ) - \eta(\varepsilon, k) ) (1 - 3 \varepsilon ) \leq 
A_k^X ( k+\beta  , S_k )
\\&{}&
 \leq 
( n_k ( k+\beta , 0 , S_k ) + \eta(\varepsilon, k)  ) (1 - 3 \varepsilon ) 
\end{eqnarray*}
since $X \subset A_k^*$ and $m(X) = 1 - 3 \varepsilon$. 

\subsection{Bounds for image areas} 
We seek an upper bound on $A_k^Y ( k+\beta  , f )$ for certain measurable subsets $Y$ of $X$. To do so, we first obtain upper bounds on  $\ell ( \Gamma'     ) = \ell ( \Gamma''     ) $, where $\ell$ denotes the spherical arc length on $\overline{  {\mathbb C}     } $. 

Let $\Gamma_j = R_k(\gamma_j ) $, for $1\leq j\leq 4$. We begin by seeking an upper bound for
$$
\ell ( \Gamma_2    ) + \ell ( \Gamma_3     ) \leq 
\ell ( \Gamma_1 \cup \Gamma_2    ) + \ell ( \Gamma_3 \cup \Gamma_4    )   .  
$$
With the ad-hoc notation 
$$
A_1 =  | T_k(   (k+\beta  ) e^{ i (  \alpha_k + \theta   ) }   )   |^{N_k-1}   | T' ( L_k(  (k+\beta  ) e^{ i (  \alpha_k + \theta   ) } )  )   |   |L_k' (  (k+\beta  ) e^{ i (  \alpha_k + \theta   ) } )     | 
$$
we have
\begin{eqnarray} \label{44a}
&{}& 
\ell ( \Gamma_3 \cup \Gamma_4    ) = \int_{\theta_3}^{\theta_5} 
\frac{  | R_k' (    (k+\beta  ) e^{ i (  \alpha_k + \theta   ) }      ) |   (  k+\beta  ) \, d\theta   } {   1 +   | R_k( (  k+\beta  ) e^{ i (  \alpha_k + \theta   ) }   )   |^{2  }         }   
 \\ &{}& 
=  N_k(  k+\beta  )  \int_{\theta_3}^{\theta_5} \frac{ A_1   \, d\theta      } {  1 +   | T_k(  ( k+\beta  ) e^{ i (  \alpha_k + \theta   ) }   )   |^{2 N_k }         }   
\notag \\ &{}& 
\leq 2 N_k(  k+\beta  )  \int_{\theta_3}^{\theta_5}  \frac{   d\theta    } {  | T_k(   (k+\beta  ) e^{ i (  \alpha_k + \theta   ) }   )   |^{N_k+1}       }   .   \notag 
\end{eqnarray}
since $|T' (z)| \leq 3$ for $|z|\leq 1 + 3 \delta$ and $|L_k'|  \equiv 2/3$.  

For $\theta_3 \leq \theta \leq \theta_5$, define $u$ in $[0,3\delta]$ by
$$
|   (k+\beta  ) e^{ i (  \alpha_k + \theta   ) }   -   k  e^{ i  \alpha_k      }     | = \frac{3}{2} (1 + u) ,
$$
or equivalently
$$
|   L_k (  (k+\beta  ) e^{ i (  \alpha_k + \theta   ) } )          | = 1+u , \qquad 0\leq u\leq 3 \delta .
$$

From the Law of Cosines we have for $\theta_3 \leq \theta \leq \theta_5$
$$
\left(    \frac{3}{2} (1 + u)     \right)^2  = g( \theta ) :=  (k+\beta  )^2 + k^2 - 2 k (k+\beta  ) \cos \theta  .
$$
Thus $  \frac{3}{2} (1 + u) = (g( \theta ))^{1/2} $ and
$$
du = \frac{1}{3} \frac{  g'( \theta )  } {   (g( \theta ))^{1/2}   } \, d\theta =  \frac{2}{3} \frac{   k (k+\beta  ) \sin \theta   } {  (g( \theta ))^{1/2}    } \, d\theta . 
$$
From (\ref{7}) we conclude that 
$$
| T_k(   (k+\beta  ) e^{ i (  \alpha_k + \theta   ) }   )   |
=
| T(L_k(   (k+\beta  ) e^{ i (  \alpha_k + \theta   ) }  ) )   | \geq 1 + \frac{u}{3} ,  
$$
implying for  $\theta_3 \leq \theta \leq \theta_5$   that 
\begin{equation} \label{45a}
 \frac{   1   } {  | T_k(   (k+\beta  ) e^{ i (  \alpha_k + \theta   ) }   )   |^{N_k+1}       }  
 \leq \frac{  1  } {  \left( 1 + \frac{u}{3} \right)^{N_k+1}       }  .  
\end{equation}
We observe that
\begin{eqnarray} \label{46a}
&{}& 
6 > 6 - 6 (1+\delta)^{-N_k} 
= 2 N_k \int_0^{3\delta} \frac{ du } {  \left( 1 + \frac{u}{3} \right)^{N_k+1}  } 
\notag \\ &{}& 
\geq  N_k \int_0^{3\delta} \frac{  \frac{3}{2} (1 + u) \, du } {  \left( 1 + \frac{u}{3} \right)^{N_k+1}  } 
\notag \\ &{}& 
\geq \frac{2}{3} N_k   \int_{\theta_3}^{\theta_5} \frac{   (g( \theta ))^{1/2}   k (k+\beta  ) \sin \theta  \, d\theta   } {    | T_k(   (k+\beta  ) e^{ i (  \alpha_k + \theta   ) }   )   |^{N_k+1}    (g( \theta ))^{1/2}    } 
\notag \\ &{}& 
\geq \frac{2}{3} N_k  k (k+\beta  ) \sin \theta_3   \int_{\theta_3}^{\theta_5} 
\frac{  d\theta     } {     | T_k(   (k+\beta  ) e^{ i (  \alpha_k + \theta   ) }   )   |^{N_k+1}        }  .  
\end{eqnarray}

Again applying the Law of Cosines, we have
\begin{eqnarray*}  
\left( \frac{3}{2} \right)^2 & = & 
k^2 +  (k+\beta  )^2  - 2 k (k+\beta  ) \cos \theta_3
\\ & = & 2 k (k+\beta  ) ( 1 - \cos \theta_3 )  + \beta^2 
\\ & = & 4  k (k+\beta  ) \sin^2 \frac { \theta_3 } {2}  + \beta^2  
\\ & \leq & 9 k^2 \sin^2 \theta_3 + \beta^2   ,
\end{eqnarray*}
yielding
$$
 \sin \theta_3  \geq \frac{  \sqrt{  \left( \frac{3}{2} \right)^2   - \beta^2     }     } {   3 k       }  .  
$$

Combining this with (\ref{46a}) we have
$$
6 > \frac{2}{9} N_k  (k+\beta  )  \sqrt{  \left( \frac{3}{2} \right)^2   - \beta^2     }  
 \int_{\theta_3}^{\theta_5} \frac{  d\theta     } {     | T_k(   (k+\beta  ) e^{ i (  \alpha_k + \theta   ) }   )   |^{N_k+1}        } ,  
$$
yielding by (\ref{36a}) and (\ref{44a}) 
\begin{equation} \label{47a}
\ell (\Gamma_3 ) \leq \frac{ 54   } {   \sqrt{  \left( \frac{3}{2} \right)^2   - \beta^2     }    } \leq \frac{  36  } {  \varepsilon^{1/4}   } .  
\end{equation}

We also have 
\begin{eqnarray} \label{48a}
&{}& 
\ell ( \Gamma_1 \cup \Gamma_2    )  
=  N_k(  k+\beta  )  \int_{\theta_1}^{\theta_3} \frac{ A_1   \, d\theta      } {  1 +   | T_k(  ( k+\beta  ) e^{ i (  \alpha_k + \theta   ) }   )   |^{2 N_k }         }   
\notag \\ &{}& 
\leq 2 N_k(  k+\beta  )  \int_{\theta_1}^{\theta_3} 
| T_k(   (k+\beta  ) e^{ i (  \alpha_k + \theta   ) }   )   |^{N_k-1}   \, d\theta  
\end{eqnarray}
as before. For $\theta_1 \leq \theta \leq \theta_3$ we define $u$ in $[0,3\delta]$ by
$$
|   (k+\beta  ) e^{ i (  \alpha_k + \theta   ) }   -   k  e^{ i  \alpha_k      }     | = \frac{3}{2} (1 - u) ,
$$
or equivalently
$$
|   L_k (  (k+\beta  ) e^{ i (  \alpha_k + \theta   ) } )          | = 1 - u , \qquad 0\leq u\leq 3 \delta .
$$

For $\theta_1 \leq \theta \leq \theta_3$ we have from the Law of Cosines 
$$
\left(    \frac{3}{2} (1 - u)     \right)^2  = g( \theta ) :=  (k+\beta  )^2 + k^2 - 2 k (k+\beta  ) \cos \theta  ,
$$
yielding $  \frac{3}{2} (1 - u) = (g( \theta ))^{1/2} $ and
$$
- du = \frac{1}{3} \frac{  g'( \theta )  } {   (g( \theta ))^{1/2}   } \, d\theta =  \frac{2}{3} \frac{   k (k+\beta  ) \sin \theta   } {  (g( \theta ))^{1/2}    } \, d\theta . 
$$
From (\ref{8}) we have 
$$
| T_k(   (k+\beta  ) e^{ i (  \alpha_k + \theta   ) }   )   |
=
| T(L_k(   (k+\beta  ) e^{ i (  \alpha_k + \theta   ) }  ) )   | \leq 1 - \frac{u}{3} . 
$$
Thus
\begin{eqnarray} \label{49a}
&{}& 
\qquad  3 > 3 - 3 (1-\delta)^{N_k} 
=  N_k \int_0^{3\delta}   \left( 1 - \frac{u}{3} \right)^{N_k-1}  \,  du 
\\ &{}& 
\geq  N_k \int_0^{3\delta} (1-u)    \left( 1 - \frac{u}{3} \right)^{N_k-1}      \,  du
\notag \\ &{}& 
\geq -  \frac{2}{3}  k (k+\beta  ) N_k   \int_{\theta_3}^{\theta_1}  \frac{ 2   (g( \theta ))^{1/2}       | T_k(   (k+\beta  ) e^{ i (  \alpha_k + \theta   ) }   )   |^{N_k-1}     \sin\theta \, d\theta } { 3 (g( \theta ))^{1/2}    } 
\notag \\ &{}& 
\geq \frac{4}{9} k (k+\beta  ) N_k   \sin \theta_1   \int_{\theta_1}^{\theta_3} 
  | T_k(   (k+\beta  ) e^{ i (  \alpha_k + \theta   ) }   )   |^{N_k-1}    \,  d\theta  . \notag 
\end{eqnarray}

From the Law of Cosines we have 
\begin{eqnarray*} 
&{}&  
\left(    \frac{3}{2} (1 - 3\delta)     \right)^2  = k^2 +  (k+\beta  )^2   - 2 k (k+\beta  ) \cos \theta_1   
\\ &{}&
\leq 9 k^2 \sin^2 \theta_1 + \beta^2 
\end{eqnarray*}
as before, implying by  (\ref{11}) and (\ref{36a})  that 
$$
 \sin \theta_1  \geq \frac{ \frac{3}{2}  \sqrt{  (1 - 3\delta)^2 - (1 - \sqrt{\varepsilon} )^2     }     } {   3 k       } \geq \frac{ \varepsilon^{1/4}    } {  2 k      }   .  
$$
Combining this with  (\ref{48a}) and (\ref{49a}) we conclude that 
$$
\ell (\Gamma_2) \leq \frac{27}{ \varepsilon^{1/4} }  . 
$$
Thus
$$
\ell (\Gamma_2 \cup \Gamma_3) \leq \frac{63}{ \varepsilon^{1/4} }  . 
$$

Since $L(z)=1/z$ is a rotation of the Riemann sphere, we see that 
$1/R_k$ maps $\gamma_2 \cup \gamma_3$ to a curve of length
$ \ell (\Gamma_2 \cup \Gamma_3) \leq  63/ \varepsilon^{1/4}    . 
$ 
Recall that $S_k(z) = 1 - 1/R_k(z)$. It is elementary that
$$
\frac{ 1 + x^2  } {  1 + (x+1)^2   } \geq \frac{  5 - \sqrt{5}  } {   5 + \sqrt{5}    } \qquad \text{ for }\, x \geq 0 .  
$$
It follows easily (recall that $\Gamma' = S_k(\gamma_2 \cup \gamma_3)$) that 
$$
\ell (\Gamma' ) \leq \frac{  5 + \sqrt{5}  } {   5 - \sqrt{5}    } \cdot \frac{63}{ \varepsilon^{1/4} }  < \frac{ 165 }{ \varepsilon^{1/4} } . 
$$
Since $\Gamma''$ is the curve conjugate to $\Gamma'$, we have
$$
\ell (\Gamma'' ) = \ell (\Gamma' )   < \frac{ 165 }{ \varepsilon^{1/4} } . 
$$

We now consider an odd $k$. Let $w_1,w_2,\dots ,w_p$ be consecutive points in $\Gamma' = S_k( \gamma_2 \cup \gamma_3  )$ with $w_1=S_k(z_2)$ and $w_p=S_k(z_4)$. By  (\ref{21}) and (\ref{23}), $ |w_1-1| = (1-\delta)^{-N_k}$ and $ |w_p-1| =  (1+\delta)^{-N_k}$. Let $\Gamma_j^*$ be the part of $\Gamma'$ from $w_j$ to $w_{j+1}$, for $1\leq j\leq p-1$. Choose $w_j$ so that $\ell (\Gamma_j^* ) = \sqrt{\varepsilon}/2$ for $1\leq j\leq p-2$ and $\ell (\Gamma_{p-1}^*) \leq \sqrt{\varepsilon}/2$. Evidently
$$
(p-2) \sqrt{\varepsilon}/2 < 165 / \varepsilon^{1/4} ,
$$
or $p-1<331/ \varepsilon^{3/4}$. 
Let $\Delta_j = D(w_j,  \sqrt{\varepsilon} )$, for $1\leq j\leq p-1$. 

Define
$$
W_a = \bigcup_{w\in \Gamma'} B(w,|w| \sqrt{\varepsilon}/2 ) 
\subset  \bigcup_{w\in \Gamma'} D(w, \sqrt{\varepsilon}/2 ) , 
$$
where we have applied (\ref{2}). Suppose that $z\in \gamma_2 \cup \gamma_3 $ with
$w = S_k(z) \in \Gamma'$. By (\ref{33}) 
$$
f(z) \in B( S_k(z) , |S_k(z) | \sqrt{\varepsilon}/2 ) \subset W_a  .
$$
Thus $ f (   \gamma_2 \cup \gamma_3  ) \subset W_a$. 

Suppose that $w'\in W_a$ with $w' \in B(w,|w| \sqrt{\varepsilon}/2 ) $ where $w\in \Gamma'$. Suppose that $w\in \Gamma_j^*$ where $1\leq j\leq p-1$. Certainly $k(w,w_j)\leq  \sqrt{\varepsilon}/2$. Since $w'\in D(w, \sqrt{\varepsilon}/2 )$, we have
$k(w',w_j)<  \sqrt{\varepsilon} $, implying that $w' \in \Delta_j$. Thus
$W_a\subset \bigcup_{j=1}^{p-1} \Delta_j\setminus \{ \infty  \}.$
Hence
$$
m(W_a) \leq \sum_{j=1}^{p-1} m( \Delta_j ) \leq (p-1) \varepsilon \leq 331 \varepsilon^{1/4} .
$$

Similarly we have
$$
f ( \tilde{\gamma}_3 \cap   \tilde{\gamma}_2      ) \subset W_b := \bigcup_{w\in \Gamma''} B(w,|w| \sqrt{\varepsilon}/2 ) 
$$
with $m(W_b) \leq 331 \varepsilon^{1/4} $. Let $\tilde{W}_k = W_a \cup W_b$. We have $m( \tilde{W}_k ) \leq 662 \varepsilon^{1/4} $.

Recall that
$$
\Gamma'''' = S_k(  \gamma_b   ) \subset \tilde{\Omega}_2
$$
by (\ref{21}). By (\ref{33}) we have
$$
f ( \gamma_b ) \subset \bigcup_{w\in \Gamma''''} B(w,|w| \sqrt{\varepsilon}/2 ) 
\subset \bigcup_{w\in \Gamma''''} D(w, \sqrt{\varepsilon}/2 ) \setminus \{ \infty \}
\subset \Omega_2  \setminus \{ \infty \}      . 
$$

Finally, $\Gamma''' = S_k(\gamma_a) \subset \tilde{\Omega}_3$ by (\ref{23}). By (\ref{33}) we conclude as above that $f(\gamma_a) \subset \Omega_3$. 

From these observations we conclude that $f(\gamma)$ is homotopic to $\Gamma = S_k(\gamma)$ in the set
$(  \tilde{W}_k  \cup \Omega_2 \cup \Omega_3     ) \setminus \{ \infty \}$, the homotopy being given by $\lambda S_k + (1-\lambda) f$, for $0\leq \lambda \leq 1$. 
Thus  $f (\gamma)$ and $\Gamma$ have the same winding number about all points in the complement of $\tilde{W}_k \cup \Omega_2 \cup \Omega_3$ and by the Argument Principle
$$
n_k ( k+\beta , a , S_k ) = n_k ( k+\beta , a , f ) 
$$
for all $a$ in ${\mathbb C} \setminus (    \tilde{W}_k  \cup \Omega_2 \cup \Omega_3       ) \supset X \setminus   \tilde{W}_k $. Trivially, even if $0 \in   \tilde{W}_k $, we have 
$$
n_k ( k+\beta , 0 , S_k ) = n_k ( k+\beta , 0 , f ) . 
$$
Combined with (\ref{43a}) we have with $r=k+\beta$
\begin{equation} \label{50a}
n_k(r,0,f)  -  \eta(\varepsilon, k)  \leq n_k(r,a,f) \leq n_k(r,0,f) + \eta(\varepsilon, k)
\end{equation}
for all $a$ in $X\setminus \tilde{W}_k $. 

If $j$ is even, then $S_j$ is a factor in the numerator of $1/f$. The above analysis shows that 
$$
n_j (j+\beta , 0, 1/f ) - \eta(\varepsilon, j)  \leq n_j ( j+\beta , a , 1/f ) \leq n_j (j+\beta , 0, 1/f ) + \eta(\varepsilon, j) 
$$
for all $a$ in $X\setminus \tilde{W}_j $  for a set $ \tilde{W}_j $ with $m( \tilde{W}_j ) < 662  \varepsilon^{1/4} $.  This can be rewritten with $r=j+\beta$ as
\begin{equation} \label{51a}
n_j(r,\infty ,f)-\eta(\varepsilon, j) \leq n_j(r,a,f) \leq n_j(r, \infty ,f)+\eta(\varepsilon, j) 
\end{equation}
for all $a$ in $X\setminus \tilde{W}_j^* $
where   $m( \tilde{W}_j^* ) < 662  \varepsilon^{1/4} $ and
$ \tilde{W}_j^* = L(\tilde{W}_j) $ with $L(z)=1/z$.   

\subsection{The final estimates}
For an odd $k\geq 7$ we consider the interval $I=I_1\cup I_2$ where
$I_1 = \left[k-\frac{1}{2} +\frac{3}{2} \sqrt{ \varepsilon } , k+\frac{1}{2} - \frac{3}{2} \sqrt{ \varepsilon } \right] $ and $I_2 = \left[k + \frac{1}{2} - \frac{3}{2} \sqrt{ \varepsilon } , k+\frac{1}{2} + \frac{3}{2} \sqrt{ \varepsilon } \right]  $. First suppose that $r\in I_1$. We have $U_{k-2} \subset \overline{B}(0,r)$. For all $a$ in $X$ we conclude from (\ref{31a}) and (\ref{32a})  for $5\leq j\leq k-2$ that 
$n_j(r,a,f) = N_j$. Also $ \overline{B}(0,r)\cap U_{k+2} = \emptyset$, implying for $j\geq k+2$ that $n_j(r,a,f)=0$ for all $a$ in $\overline{{\mathbb C}}$ and thus $A_j (r,f)=0$.

Recall that $f$ has a pole of multiplicity $N_k$ at the special point $ (k-\frac{1}{2}) e^{ i  \alpha_k } $ in $D_k$, with $N_k$ simple zeros on $\partial D_k$. For $j=k-1$ or $j=k+1$, the function $f$ has a zero of multiplicity $N_j$ at the special point in $D_j$ with $N_j$ simple poles on $\partial D_j$. We note that for $k-1\leq j\leq k+1$, $r=j+\beta$ for some $\beta$ satisfying (\ref{36a}). Thus (\ref{50a}) is applicable for our $r$ and (\ref{51a}) is applicable for our $r$ with $j=k-1$ or $j=k+1$. 

Let 
$$
X^* = X \setminus (  \tilde{W}_{k-1}^*  \cup \tilde{W}_k  \cup \tilde{W}_{k+1}^*   ) 
$$
where the $\tilde{W}$ are the exceptional sets associated with $r$ as in (\ref{50a}) and (\ref{51a}). We have
\begin{equation} \label{52a}
1 - 3 \varepsilon -  1986  \varepsilon^{1/4}   \leq m(X^*) \leq 1 - 3 \varepsilon  .  
\end{equation}

From (\ref{34a}), (\ref{35a}), (\ref{50a}), (\ref{51a}), and (\ref{52a}),    we have
\begin{eqnarray*} 
&{}&
A(r,f)  \leq  A^{X^*} (r,f) + n(r) (1 - m(X^*) ) 
\\ &{}&
= 
\sum_{j=5}^{k+1} A_j^{X^*} (r,f) + n(r) (1 - m(X^*) ) 
\\ &{}&
\leq \left(  \sum_{j=5}^{k-2}  N_j  + n_{k-1}(r,\infty,f) + n_{k}(r,0,f) + n_{k+1}(r,\infty,f)      + \eta^*(\varepsilon, k)  \right)  m(X^*)
\\ &{}&
  + n(r) ( 3 \varepsilon +  1986  \varepsilon^{1/4} )    , 
\end{eqnarray*} 
where $\eta^*(\varepsilon, k) = \eta (\varepsilon, k-1) + \eta(\varepsilon, k) + \eta (\varepsilon, k + 1)$. 
Thus
\begin{eqnarray}  \label{53a}
&{}&
\frac{  n(r,\infty,f)   } {   A(r,f) - ( \eta^*(\varepsilon, k) + n(r) (    3 \varepsilon +  1986  \varepsilon^{1/4}   )           )      } 
\\ &{}&
  \geq \frac{ \sum_{j=5}^{k-2}  N_j  + n_{k-1}(r,\infty,f) + N_k + n_{k+1}(r,\infty,f)      } {    \sum_{j=5}^{k-2}  N_j  + n_{k-1}(r,\infty,f) + n_{k}(r,0,f) + n_{k+1}(r,\infty,f)      }    . 
  \notag  
\end{eqnarray} 

Since the right hand side of (\ref{53a}) is larger than $1$, it is decreased by replacing any of $n_{k-1}(r,\infty,f)$, $n_{k}(r,0,f)$,  or $n_{k+1}(r,\infty,f) $ by something larger.

Note that the disk $B(0,r)$ is tangent to the line $ {\rm Re}\, (  z e^{ -i \alpha_j }    ) =r$ at $z= r e^{ i \alpha_j } $ and that $ \overline{B}(0,r)$ lies in the half plane
$$
H_j^* := \{ z \colon  {\rm Re}\, (  z e^{ -i \alpha_j }    ) \leq  r \}  .
$$
Let $n_j^* (r,a,f)$ denote the number of solutions of $f(z)=a$ in $H_j^* \cap U_j$ (rather than in $\overline{B}(0,r) \cap U_j $) counting multiplicities. Thus $n_j^*$ is (slightly) larger than $n_j$. From elementary geometry it follows that
\begin{equation} \label{54a}
n_j(r, a ,f) \geq n_j^* \left(r - \frac{3}{2} \sqrt{ \varepsilon } ,a,f \right)  
\end{equation}
for all $a$ if $r>r_0(\varepsilon)$. Evidently we may replace $n$ by $n^*$ and $r$ by $ k+\frac{1}{2} $ on the right hand side of (\ref{53a}), obtaining 
\begin{eqnarray}  \label{55a}
&{}& \qquad  \qquad 
\frac{  n(r,\infty,f)   } {   A(r,f) - ( \eta^*(\varepsilon, k) + n(r) (    3 \varepsilon +  1986  \varepsilon^{1/4}   )           )      } 
\\ &{}&
  \geq \frac{ \sum_{j=5}^{k-2}  N_j  + n_{k-1}^* ( k+\frac{1}{2} ,\infty,f) + N_k + n_{k+1}^* ( k+\frac{1}{2} ,\infty,f)      } {    \sum_{j=5}^{k-2}  N_j  + n_{k-1}^* \left( k+\frac{1}{2} ,\infty, f \right) + n_{k}^* \left( k+\frac{1}{2} ,0, f \right) + n_{k+1}^* \left( k+\frac{1}{2} ,\infty, f \right)      }    . 
  \notag  
\end{eqnarray} 

Now $L_j$ maps the line $  {\rm Re}\, (  z e^{ -i \alpha_j }    ) = k+\frac{1}{2} $ to the line $  {\rm Re}\, w = \frac{2}{3} (k-j) + \frac{1}{3} $. In particular, $L_{k-1}$ maps the line $  {\rm Re}\, (  z e^{ -i \alpha_{k-1} }    ) = k+ \frac{1}{2} $ to the line $  {\rm Re}\, w = 1 $; $L_k$ maps the line $  {\rm Re}\, (  z e^{ -i \alpha_k }    ) = k+\frac{1}{2} $ to the line $  {\rm Re}\, w =  \frac{1}{3} $; and $L_{k+1}$ maps the line $  {\rm Re}\, (  z e^{ -i \alpha_{k+1} }    ) = k+\frac{1}{2} $ to the line $  {\rm Re}\, w = - \frac{1}{3} $. 

Trivially we have
\begin{equation} \label{56a}
 n_{k-1}^* \left(k+\frac{1}{2}   , \infty , f \right)  = N_{k-1}  . 
\end{equation}
From (\ref{9a}), (\ref{11a}), and the observations in the paragraph above we have for large $k$
\begin{equation} \label{57a}
 n_{k}^* \left(k+\frac{1}{2}   , 0 , f \right)  \leq 0.78366 \,  N_{k}  . 
\end{equation}
Likewise from (\ref{8a}) and (\ref{10a}) we conclude that 
\begin{equation} \label{58a}
 n_{k+1}^* \left(k+\frac{1}{2}   , \infty , f \right)  \leq 0.60818 \,  N_{k+1}  . 
\end{equation}

Increasing $N_j = \left[  C^{ k_0 + j      }        \right] $ to $C^{ k_0 + j      }  $, then canceling a common factor of $C^{k_0}$ in the numerator and denominator, and noting  the change in the fifth decimal place in the coefficient of $C^k$ in the denominator to account for the change from $N_k$ to $C^{k_0 +k}$ in the numerator, we have from (\ref{55a}), (\ref{56a}), (\ref{57a}), and (\ref{58a}), 
\begin{eqnarray}  \label{59a}
&{}& \qquad  \qquad 
\frac{  n(r,\infty,f)   } {   A(r,f) - ( \eta^*(\varepsilon, k) + n(r) (    3 \varepsilon +  1986  \varepsilon^{1/4}   )           )      } 
\\ &{}&
  \geq \frac{ \frac{ C^{k-1}} { C-1 } - \frac{ C^{5}} { C -1 }  + C^{k-1} + C^k + 0.60818 C^{k+1}     } {   \frac{ C^{k-1}} { C -1 } - \frac{ C^{5}} { C -1 }  + C^{k-1}   + 0.78367 C^k  +  0.60818 C^{k+1}    }    
  \notag  
\\ &{}&
\geq 1 + \frac{  0.21633   } {  \frac{ 1 } { C -1 }   + 0.78367 + 0.60818 C          }  , \notag  
\end{eqnarray} 
where in the last step we have added $C^5/(C -1)$ to both numerator and denominator.

It is elementary that the last expression is maximized at $C=1+( 0.60818 )^{-1/2} \approx 2.28228 .$ With $C = 2.28228$, this expression is greater than $1.07329$. 
With this choice of $C$, we obtain for large $r\in I_1$ that
$$
\frac{  n(r )   } {   A(r,f) - ( \eta^*(\varepsilon, k) + n(r) (    3 \varepsilon +  1986  \varepsilon^{1/4}   )           )      }  > 1.07329  ,  
$$
which we rearrange to
\begin{equation} \label{60a}
\frac{  n(r )   } {   A(r,f) } > 
\frac{ 1.07329 - \eta^*(\varepsilon, k) (1.07329 ) /A(r,f)        } {  1 + 1.07329 (    3 \varepsilon +  1986  \varepsilon^{1/4}   )        } > 1.07328  , 
\end{equation} 
valid for $r$ in $I_1 = \left[k-\frac{1}{2} +\frac{3}{2} \sqrt{ \varepsilon } , k+\frac{1}{2} - \frac{3}{2} \sqrt{ \varepsilon } \right] $ for large odd $k$ and small $  \varepsilon >0$. 
Note that in the last step we have used $A(r,f) > A^X_{k-2} (r,f)  = N_{k-2} (1- 3 \varepsilon)$.  

We now suppose that $r \in I_2$. We note that $ \overline{B}(0,r) \cap U_{k+3} = \emptyset$ and thus $A_j(r,f)=0$ for $j\geq k+3$. 
For $k\leq j\leq k+2$, $L_j$ maps the lines $  {\rm Re}\, (  z e^{ -i \alpha_{j} }    ) = k+\frac{1}{2} + \frac{3}{2} \sqrt{ \varepsilon } $ to the lines $  {\rm Re}\, w =  \frac{1}{3} +  \sqrt{ \varepsilon } $, $  {\rm Re}\, w = -  \frac{1}{3} +  \sqrt{ \varepsilon } $, 
and $  {\rm Re}\, w = -  1 +  \sqrt{ \varepsilon } $, respectively. From (\ref{12a}) with $s= \sqrt{ \varepsilon } $ and from (\ref{57a}) and (\ref{58a})  
we have for large $k$ 
\begin{equation} \label{62a}
n_{k}^* \left( k+\frac{1}{2} + \frac{3}{2} \sqrt{ \varepsilon } , 0 , f \right) < \left( 0.78366 + \frac{4} {\pi} \varepsilon^{1/4}   \right) N_{k} ,  
\end{equation}
\begin{equation} \label{63a}
n_{k+1}^* \left( k+\frac{1}{2} + \frac{3}{2} \sqrt{ \varepsilon } , \infty , f \right)  < \left( 0.60818 + \frac{4} {\pi} \varepsilon^{1/4}    \right) N_{k+1} , 
\end{equation}
and 
\begin{equation} \label{64a}
n_{k+2}^* \left( k+\frac{1}{2} + \frac{3}{2} \sqrt{ \varepsilon } , 0 , f \right)  <  \frac{4} {\pi} \varepsilon^{1/4}  N_{k+2} .  
\end{equation}

For $k\leq j\leq k+2$ we note that $r' = k+\frac{1}{2} + \frac{3}{2} \sqrt{ \varepsilon } $ is of the form $r' = j+\beta$ for some $\beta$ satisfying (\ref{36a}). Thus (\ref{50a}) and (\ref{51a}) are valid with $r' = k+\frac{1}{2} + \frac{3}{2} \sqrt{ \varepsilon } $ for $a$ in a set $X^{**} = X \setminus W^{**}$ where $W^{**}$ is a union of exceptional sets associated with $S_j$, for $k\leq j\leq k+2$, on $S(0,r')$ and $m(W^{**}) < 1986 \varepsilon^{1/4} $. 

We have 
\begin{eqnarray*} 
&{}& 
A(r,f) \leq A^{ X^{**} } (r,f) + n(r) ( 1 - m(  X^{**}   ) ) 
\\ &{}&
\leq A^{ X^{**} } (r' , f ) + n(r) ( 1 - m(  X^{**}   ) )
\\ &{}&
= \sum_{j=5}^{k+2} A_j^{ X^{**} } \left(  k+\frac{1}{2} + \frac{3}{2} \sqrt{ \varepsilon } , f \right) 
+ n(r) ( 1 - m(  X^{**}   ) ) .  
\end{eqnarray*} 
Since $U_{k-1} \subset B\left(0, k+\frac{1}{2} + \frac{3}{2} \sqrt{ \varepsilon } \right)$, upon combining (\ref{50a}) and (\ref{51a}) with (\ref{62a}), (\ref{63a}), and (\ref{64a}) we obtain
\begin{eqnarray}  \label{65a}
&{}& \qquad  \qquad \qquad 
A(r,f) \leq \Biggl\{ \sum_{j=5}^{k-1} N_j + \left( 0.78366 +  \frac{4} {\pi} \varepsilon^{1/4} \right) N_k 
\\ &{}&
+ \left( 0.60818 +  \frac{4} {\pi} \varepsilon^{1/4}  \right) N_{k+1} + 
 \frac{4} {\pi} \varepsilon^{1/4}  N_{k+2} +  \eta^{**}(\varepsilon, k) \Biggr\} m(X^{**})  \notag
\\ &{}&
 + n(r) (  3 \varepsilon +  1986  \varepsilon^{1/4}   )    ,      
  \notag
\end{eqnarray}
where $\eta^{**}(\varepsilon, k) =  \eta(\varepsilon,k) + \eta (\varepsilon, k +1) + \eta ( \varepsilon, k+2) $. 

For $k-1\leq j\leq k+1$, $L_j$ maps the lines 
$  {\rm Re}\, (  z e^{ -i \alpha_{j} }    ) = k+\frac{1}{2} - \frac{3}{2} \sqrt{ \varepsilon } $ to the lines $  {\rm Re}\, w =  1 - \sqrt{ \varepsilon } $, $  {\rm Re}\, w =  \frac{1}{3} - \sqrt{ \varepsilon } $, and $  {\rm Re}\, w =  -  \frac{1}{3} - \sqrt{ \varepsilon } $, respectively. From (\ref{12a}) with $s= 3 \sqrt{ \varepsilon } $, (\ref{54a}), (\ref{56a}), (\ref{8a}), and (\ref{10a}) we have 
\begin{eqnarray} \label{66a}
&{}&
n_{k-1} \left( k+\frac{1}{2} - \frac{3}{2} \sqrt{ \varepsilon } , \infty , f \right)  \geq 
n_{k-1}^* \left( k+\frac{1}{2} - 3\sqrt{ \varepsilon } , \infty , f \right) 
\\ &{}&
> \left(1 -  \frac{8} {\pi} \varepsilon^{1/4}  \right) N_{k-1} \notag
\end{eqnarray}
and
\begin{eqnarray} \label{67a}
&{}&
n_{k+1} \left( k+\frac{1}{2} - \frac{3}{2} \sqrt{ \varepsilon } , \infty , f \right)  \geq 
n_{k+1}^* \left( k+\frac{1}{2} - 3 \sqrt{ \varepsilon } , \infty , f \right)
 \\ &{}&
> \left(0.60817 -  \frac{8} {\pi} \varepsilon^{1/4}  \right) N_{k+1} .  \notag
\end{eqnarray}
Thus for $r\in I_2$ with large $k$, from (\ref{65a}), (\ref{66a}), and (\ref{67a}) we have 
\begin{eqnarray}  \label{68a}
&{}& \qquad  \qquad \qquad 
\frac{  n(r )   } {   A(r,f) - ( \eta^{**}(\varepsilon, k) + n(r) (    3 \varepsilon +  1986  \varepsilon^{1/4}   )           )      } 
\\ &{}&
 \geq \frac{   n \left( k+\frac{1}{2} - \frac{3}{2} \sqrt{ \varepsilon } , \infty , f \right)   } {     A(r,f) - ( \eta^{**}(\varepsilon, k) + n(r) (    3 \varepsilon +  1986  \varepsilon^{1/4}   )           )           } \notag 
\\ &{}&
\geq \frac{ \sum_{j=5}^{k-2} N_j + \left( 1 - \frac{8} {\pi} \varepsilon^{1/4} \right) N_{k-1}
+ N_k + \left(  0.60817 -  \frac{8} {\pi} \varepsilon^{1/4}      \right) N_{k+1}  
    } {     \sum_{j=5}^{k-1} N_j + \left( 0.78366 + \frac{4} {\pi} \varepsilon^{1/4} \right) N_{k}
+  \left(  0.60818 + \frac{4} {\pi} \varepsilon^{1/4}      \right) N_{k+1} +   \frac{4} {\pi} \varepsilon^{1/4}       N_{k+2}      }          
.  \notag
\end{eqnarray} 

With the choice $C=2.28228$, our previous discussion (see (\ref{59a})) shows that for small $\varepsilon >0$ the right hand side of (\ref{68a}) is greater than $1.073284$. As before we rearrange  (\ref{68a})  to obtain  for small $\varepsilon >0$ 
\begin{equation} \label{69a}
\frac{  n(r )   } {   A(r,f) } > 1.07328 , \quad r\in I_2, \quad k \text{ odd and large} .  
\end{equation}

Combining (\ref{60a}) and (\ref{69a}) we have
\begin{equation} \label{70a}
\frac{  n(r )   } {   A(r,f) } > 1.07328 , \quad r\in I = \left[ k-\frac{1}{2}+ \frac{3}{2} \sqrt{ \varepsilon }, k+\frac{1}{2} + \frac{3}{2} \sqrt{ \varepsilon }   \right]  
\end{equation}
where $k$ is odd and large. 

Since $k+1$ is even and $S_{k+1}$ is a factor in the numerator of $1/f$, we apply the above analysis to $1/f$ to conclude with obvious notation
$$
\frac{  n(r,1/f )   } {   A(r,1/f) } > 1.07328 , \quad r\in I = \left[ k+\frac{1}{2}+ \frac{3}{2} \sqrt{ \varepsilon }, k+\frac{3}{2} + \frac{3}{2} \sqrt{ \varepsilon }   \right]  .  
$$
Since $n(r,1/f)=n(r)$ and $A(r,1/f)=A(r,f)$ we now have
$$
\frac{  n(r )   } {   A(r,f) } > 1.07328 , \quad   k-\frac{1}{2}+ \frac{3}{2} \sqrt{ \varepsilon } \leq r \leq  k+\frac{3}{2} + \frac{3}{2} \sqrt{ \varepsilon }   ,  
$$
for all large odd $k$  and some fixed $\varepsilon > 0$. We conclude that
\hfil\break
 $\liminf_{r\to \infty } \frac{n(r)}{A(r,f)} \geq 1.07328$. This completes the proof of Theorem~\ref{th1}.

\end{document}